\documentclass[runningheads,a4paper]{llncs}

\usepackage{amssymb}
\usepackage{graphicx}

\usepackage{url} 
\newcommand{\keywords}[1]{\par\addvspace\baselineskip
\noindent\keywordname\enspace\ignorespaces#1}

\begin{document}

\mainmatter  %

\title{Computing the Chow Variety   \\ of Quadratic Space Curves}

\titlerunning{Computing the Chow Variety of Quadratic Space Curves}

\author{Peter B\"urgisser\inst{1}%
\thanks{This author was partially supported by DFG grant BU 1371/2-2.}
\and Kathl\'en Kohn\inst{1}%
\thanks{This author was supported by a Fellowship from the Einstein Foundation Berlin.}
\and Pierre Lairez\inst{1}\textsuperscript{$\star$}\and Bernd Sturmfels\inst{1,2}\textsuperscript{$\star\star$}}
\authorrunning{Computing the Chow Variety of Quadratic Space Curves}

\institute{Institute of Mathematics, Technische Universit\"at Berlin, Germany
\and Department of Mathematics, University of California, Berkeley, USA}

\toctitle{Computing the Chow Variety of Quadratic Space Curves}
\tocauthor{Peter B\"urgisser, Kathl\'en Kohn, Pierre Lairez, Bernd Sturmfels}
\maketitle

\begin{abstract}
Quadrics in the Grassmannian of lines in $3$-space
form a $19$-dimensional projective space. We study the subvariety of 
coisotropic hypersurfaces. Following Gel'fand, Kapranov 
and Zelevinsky, it decomposes into Chow forms of plane conics,
Chow forms of pairs of lines, and Hurwitz forms of quadric surfaces.
We compute the ideals of these loci.
\end{abstract}

\section*{Introduction}

The Chow variety, introduced in 1937 by Chow and van der Waerden~\cite{CvW},
parameterizes algebraic cycles of any fixed dimension and degree in a
projective space,
each given by its Chow form.
 The case of curves in~$\bbbp^3$ goes back to an 1848 paper by
Cayley~\cite{cayley}.  A fundamental problem, addressed by Green and
Morrison~\cite{green_morrison} as well as Gel'fand, Kapranov and
Zelevinsky~\cite[\S4.3]{gkz}, is to describe the equations defining Chow
varieties.  We present a definitive computational solution for the smallest
non-trivial case, namely for cycles of dimension~$1$ and degree~$2$
in~$\bbbp^3$.

The \emph{Chow form} of a  cycle of degree $2$ is a quadratic form in the
Pl\"ucker coordinates of the Grassmannian~$G(2,4)$ of lines
in~$\bbbp^3$. Such a quadric in~$G(2,4)$  represents the set
 of all lines that intersect the given cycle. Quadratic forms in Pl\"ucker
coordinates form a projective space~$\bbbp^{19}$. The Chow variety we are
interested in, denoted~$G(2,2,4)$, is the set of all Chow forms in
that~$\bbbp^{19}$.  The aim of this note is to make  the concepts
in~\cite{cayley,CvW,green_morrison} and \cite[\S 4.3]{gkz} completely explicit.

We start with the~$9$-dimensional subvariety of~$\bbbp^{19}$ whose
points are the {\em coisotropic quadrics} in~$G(2,4)$.  By \cite[\S4.3, Theorem~3.14]{gkz},
this  decomposes as the Chow variety and the variety of Hurwitz
forms~\cite{Stu}, representing lines that are tangent to a 
quadric surface in~$\bbbp^3$.  Section~\ref{sec:coiso} 
studies the ideal generated by the coisotropy conditions.   
We work in a polynomial ring in~$20$ variables, 
one for each quadratic Pl\"ucker monomial on $G(2,4)$
minus one for the Pl\"ucker relation. 
We derive the coisotropic ideal from the differential characterization of
coisotropy. Proposition~\ref{prop:coisotropic} exhibits the decomposition
of this ideal into three minimal primes.
In particular, this shows that the coisotropic ideal is radical, 
and  it hence resolves the degree $2$ case of
 a problem posed in 1986 by Green and Morrison \cite{green_morrison}.
 They wrote:  `We do not know whether [the differential characterization of
coisotropy] generates the full ideal of these Chow variables.' 

Section~\ref{sec:chowvar} derives the radical ideal of
the Chow variety $G(2,2,4)$ in~$\bbbp^{19}$.
Its two minimal primes represent Chow forms of plane conics and Chow forms of pairs of
lines. We also study the characterization of Chow forms among all coisotropic
quadrics by the vanishing of certain differential forms.
These represent the integrability of the $\alpha$-distribution
in \cite[\S4.3, Theorem 3.22]{gkz}.
After saturation by the irrelevant ideal, the integrability ideal is found to be radical. 

\section{Coisotropic Quadrics}
\label{sec:coiso}

The Grassmannian $G(2,4)$ is a quadric in $\bbbp^5$.
Its points are lines in $\bbbp^3$. We represent these lines using 
dual Pl\"ucker coordinates $\vec{p} = (p_{01}, p_{02}, p_{03}, p_{12}, p_{13}, p_{23})$
subject to the Pl\"ucker relation~$p_{01}p_{23}-p_{02}p_{13} + p_{03}p_{12}$.
Following \cite[\S 2]{Stu}, by
{\em dual} coordinates we mean that
$p_{ij}$ is the $ij$-minor of a $2 \times 4$-matrix
whose rows span the line.
The generic quadric in $G(2,4)$   %
is written as a generic quadratic form
\begin{equation}
\label{general_quadric}
Q(\vec{p}) \ =\ \,  \vec{p} \cdot \left(
\begin{array}{cccccc}
c_0 & c_1 & c_2 & c_3 & c_4 & c_5 \\
c_1 & c_6 & c_7 & c_8 & c_9 & c_{10} \\
c_2 & c_7 & c_{11} & c_{12} & c_{13} & c_{14} \\
c_3 & c_8 & c_{12} & c_{15} & c_{16} & c_{17} \\
c_4 & c_9 & c_{13} & c_{16} & c_{18} & c_{19} \\
c_5 & c_{10} & c_{14} & c_{17} & c_{19} & c_{20}
\end{array} \right) \cdot \vec{p}^T \; .
\end{equation}
The quadric
$Q(\vec p)$ is an element in~$V := \bbbc [\vec{p}]_2 / \bbbc  \lbrace 
p_{01}p_{23}-p_{02}p_{13} + p_{03}p_{12} \rbrace \simeq \bbbc^{21} / \bbbc$.
Hence, $\vec{c} =  (c_0, c_1, \ldots, c_{20})$ serves as homogeneous coordinates on $\bbbp^{19} = \bbbp(V)$, which~-- due to the Pl\"ucker relation~-- need to be understood modulo
\begin{equation}\label{eq:coordinate-indet}
c_5 \mapsto c_5 + \lambda, \,\,
c_9 \mapsto c_9 - \lambda, \,\,
c_{12} \mapsto c_{12} + \lambda \; .
\end{equation}
The coordinate ring $\bbbq[V]$ is a subring of~$\bbbq[c_0,c_1,\ldots,c_{20}]$,
namely it is the invariant ring of the additive group action~(\ref{eq:coordinate-indet}).
Hence $\bbbq[V]$ is the polynomial ring in  $20$ variables
$\,c_0,c_1,c_2,c_3,c_4,c_5-c_{12},c_6,c_7,c_8, c_9+c_{12},c_{10},c_{11},c_{13},\ldots,c_{20}  $.

We are interested in the $\vec{c}$'s that lead to \emph{coisotropic}
hypersurfaces of~$G(2,4)$. For these, the tangent space
at any point~$\ell$, considered as a subspace of~$\,T_\ell G(2,4) $ $ =
\mathrm{Hom}(\ell, \bbbc^4/\ell)$, has the form $\left\{ \varphi  \
\middle|\  \varphi(a) = 0  \right\} + \left\{\varphi \ \middle|\
\mathrm{im}(\varphi) \subset M  \right\}$, for some~$a\in\ell \backslash \{0\}$
 and some plane~$M$ in~$\bbbc^4/\ell$.
By~\cite[\S 4.3, (3.24)]{gkz}, the quadric hypersurface $\{Q(\vec p) = 0\}$
in~$G(2,4)$  is coisotropic if and only if there exist~$s,t \in \bbbc$ such that
\begin{equation}\label{eq:coisotropy_differential}
\frac{\partial Q}{\partial p_{01}} \cdot \frac{\partial Q}{\partial p_{23}} - \frac{\partial Q}{\partial p_{02}} \cdot \frac{\partial Q}{\partial p_{13}} + \frac{\partial Q}{\partial p_{03}} \cdot \frac{\partial Q}{\partial p_{12}}
= s \cdot Q + t \cdot \left( p_{01}p_{23}-p_{02}p_{13} + p_{03}p_{12} \right) \; .
\end{equation}
Equivalently, the vector $(t,s,-1)^T$ is in the  kernel of the $21 {\times}
3$ matrix in Figure~\ref{fig:bigmat}. The $3\times 3$ minors of this matrix are
all in the subring~$\bbbq[V]$.  The {\em coisotropic ideal}~$I$ is the ideal
of~$\bbbq[V]$ generated by these minors.  The  subscheme~$V(I)$
of~$\bbbp^{19}=\bbbp(V)$ represents all coisotropic hypersurfaces~$\lbrace Q =
0 \rbrace$ of degree two in~$G(2,4)$.
Using computations with {Maple} and {Macaulay2}~\cite{M2}, we found that~$I$
has codimension~10, degree~92 and is minimally generated by~175 cubics.
Besides,~$V(I)$ is the reduced union of three  components, of dimensions nine,
eight and five. 
\looseness=-1

\begin{proposition} \label{prop:coisotropic}
The coisotropic ideal is  the intersection of three prime ideals:
\begin{equation}
\label{coisotropy_ideal_decomposition}
I \,\,= \,\,P_{\mathrm{Hurwitz}} \,\cap \, P_{\mathrm{ChowLines}} \,\cap \,P_{\mathrm{Squares}} \; .
\end{equation}
So, $I$ is radical.  The prime $P_{\mathrm{Hurwitz}}$ has codimension 10 and
degree~92, it is minimally generated by~20 quadrics, and its variety $V (
P_{\mathrm{Hurwitz}} )$ consists of Hurwitz forms of quadric surfaces in
$\bbbp^3$. The prime $P_{\mathrm{ChowLines}}$ has codimension~11 and degree~140,
it is minimally generated by~265 cubics, and $V (P_{\mathrm{ChowLines}} )$
consists of Chow forms of pairs of lines in $\bbbp^3$. The prime
$P_{\mathrm{Squares}}$ has codimension~14 and degree~32, it is minimally
generated by~84 quadrics, and $V ( P_{\mathrm{Squares}})$ consists of all
quadrics 
 $Q(\vec{p})$ that are squares modulo the Pl\"ucker relation. 
\end{proposition}

\begin{figure}
\centering
\[
 \left( \begin{array}{ccc}
    0 & {c}_{0} & 2  {c}_{0}  {c}_{5}-2  {c}_{1}  {c}_{4}+2 {c}_{2}  {c}_{3}\\
    0 & {c}_{1} & {c}_{0}  {c}_{10}-{c}_{1} {c}_{9}+{c}_{2}  {c}_{8}+{c}_{3}  {c}_{7}-{c}_{4} {c}_{6}+{c}_{1}  {c}_{5}\\
    0 & {c}_{2} & {c}_{0} {c}_{14}-{c}_{1}  {c}_{13}+{c}_{2}  {c}_{12}+{c}_{3} {c}_{11}-{c}_{4}  {c}_{7}+{c}_{2}  {c}_{5}\\
    0 & {c}_{3} & {c}_{0}  {c}_{17}-{c}_{1}  {c}_{16}+{c}_{2}  {c}_{15}+{c}_{3} {c}_{12}-{c}_{4}  {c}_{8}+{c}_{3}  {c}_{5}\\
    0 & {c}_{4} & {c}_{0}  {c}_{19}-{c}_{1}  {c}_{18}+{c}_{2}  {c}_{16}+{c}_{3} {c}_{13}-{c}_{4}  {c}_{9}+{c}_{4}  {c}_{5}\\
    1 & {c}_{5} & {c}_{0}  {c}_{20}-{c}_{1}  {c}_{19}+{c}_{2}  {c}_{17}+{c}_{3} {c}_{14}-{c}_{4}  {c}_{10}+{c}_{5}^{2}\\ 0 & {c}_{6} & 2 {c}_{1}  {c}_{10}-2  {c}_{6}  {c}_{9}+2  {c}_{7}  {c}_{8}\\
    0 & {c}_{7} & {c}_{1}  {c}_{14}-{c}_{6}  {c}_{13}+{c}_{7} {c}_{12}+{c}_{8}  {c}_{11}+{c}_{2}  {c}_{10}-{c}_{7} {c}_{9}\\
    0 & {c}_{8} & {c}_{1}  {c}_{17}-{c}_{6} {c}_{16}+{c}_{7}  {c}_{15}+{c}_{8}  {c}_{12}+{c}_{3} {c}_{10}-{c}_{8}  {c}_{9}\\
    -1 & {c}_{9} & {c}_{1} {c}_{19}-{c}_{6}  {c}_{18}+{c}_{7}  {c}_{16}+{c}_{8} {c}_{13}+{c}_{4}  {c}_{10}-{c}_{9}^{2}\\
    0 & {c}_{10} & {c}_{1}  {c}_{20}-{c}_{6}  {c}_{19}+{c}_{7}  {c}_{17}+{c}_{8} {c}_{14}-{c}_{9}  {c}_{10}+{c}_{5}  {c}_{10}\\
    0 & {c}_{11} & 2  {c}_{2}  {c}_{14}-2  {c}_{7}  {c}_{13}+2  {c}_{11} {c}_{12}\\
    1 & {c}_{12} & {c}_{2}  {c}_{17}-{c}_{7} {c}_{16}+{c}_{11}  {c}_{15}+{c}_{3}  {c}_{14}-{c}_{8} {c}_{13}+{c}_{12}^{2}\\
    0 & {c}_{13} & {c}_{2} {c}_{19}-{c}_{7}  {c}_{18}+{c}_{11}  {c}_{16}+{c}_{4} {c}_{14}+{c}_{12}  {c}_{13}-{c}_{9}  {c}_{13}\\
    0 & {c}_{14} & {c}_{2}  {c}_{20}-{c}_{7}  {c}_{19}+{c}_{11} {c}_{17}+{c}_{12}  {c}_{14}+{c}_{5}  {c}_{14}-{c}_{10} {c}_{13}\\
    0 & {c}_{15} & 2  {c}_{3}  {c}_{17}-2  {c}_{8} {c}_{16}+2  {c}_{12}  {c}_{15}\\
    0 & {c}_{16} & {c}_{3} {c}_{19}-{c}_{8}  {c}_{18}+{c}_{4}  {c}_{17}+{c}_{12} {c}_{16}-{c}_{9}  {c}_{16}+{c}_{13}  {c}_{15}\\
    0 & {c}_{17} & {c}_{3}  {c}_{20}-{c}_{8}  {c}_{19}+{c}_{12} {c}_{17}+{c}_{5}  {c}_{17}-{c}_{10}  {c}_{16}+{c}_{14} {c}_{15}\\
    0 & {c}_{18} & 2  {c}_{4}  {c}_{19}-2  {c}_{9} {c}_{18}+2  {c}_{13}  {c}_{16}\\
    0 & {c}_{19} & {c}_{4} {c}_{20}-{c}_{9}  {c}_{19}+{c}_{5}  {c}_{19}-{c}_{10} {c}_{18}+{c}_{13}  {c}_{17}+{c}_{14}  {c}_{16}\\
    0 & {c}_{20} & 2  {c}_{5}  {c}_{20}-2  {c}_{10}  {c}_{19}+2  {c}_{14} {c}_{17}
  \end{array} \right)
  \]
  \label{fig:bigmat}
  \caption{This matrix has rank  $\leq 2$ if and only 
  if the quadric given by~$\vec c$ is coisotropic.}
\end{figure}

This proposition answers a question due to Green and Morrison, who had asked in 
 \cite{green_morrison} whether $I$ is radical.
 To derive the prime decomposition~(\ref{coisotropy_ideal_decomposition}), we computed the three prime ideals as kernels of homomorphisms of polynomial rings, each expressing the relevant
 geometric condition. This construction ensures that the ideals are prime. We then
verified that their intersection equals $I$. For details, check our computations,
 using the link given at the end of this article.

From the geometric perspective of \cite{gkz}, the third 
prime $P_{\mathrm{Squares}}$ is extraneous,
because nonreduced hypersurfaces in $G(2,4)$ are excluded 
by Gel'fand, Kapranov 
and Zelevinsky.
 Theorem 3.14 in \cite[\S4.3]{gkz} concerns irreducible hypersurfaces, and the identification of Chow forms within the coisotropic hypersurfaces~\cite[\S4.3, Theorem 3.22]{gkz} assumes
   the corresponding polynomial to be squarefree.
With this, the following would be the correct ideal for the coisotropic variety in $\bbbp^{19}$:
\begin{equation}
\label{correct_coisotropic_ideal}
P_{\mathrm{Hurwitz}} \,\cap \,P_{\mathrm{ChowLines}} 
\,\,= \,\, \left( I : P_{\mathrm{Squares}} \right) \; .
\end{equation}
This means that the reduced coisotropic quadrics in~$G(2,4)$ are either Chow forms of curves or Hurwitz forms of surfaces.
The ideal in~(\ref{correct_coisotropic_ideal}) has codimension~10, degree~92,
and is minimally generated by~175 cubics and~20 quartics in $\bbbq[V]$.

A slightly different point of view on the coisotropic ideal is presented in a recent
paper of Catanese \cite{catanese}.  He derives a variety in $\mathbb{P}^{20} = 
\mathbb{P}(\mathbb{C}[{\bf p}]_2)$ which projects isomorphically onto our variety $V(I) \subset \mathbb{P}^{19}$.
The center of projection is the Pl\"ucker quadric.
To be precise,  Proposition 4.1 in \cite{catanese} states the following:
   For every $Q \in \bbbc[\vec{p}]_2 \backslash
   \bbbc  \left( p_{01}p_{23}-p_{02}p_{13} + p_{03}p_{12} \right) 
$    satisfying~(\ref{eq:coisotropy_differential}) there is a unique $\lambda \in \mathbb{C}$ such that
the quadric $Q_\lambda := Q + \lambda \cdot \left( p_{01}p_{23}-p_{02}p_{13} + p_{03}p_{12} \right)$
satisfies
\begin{equation}\label{eq:coisotropy_catanese}
\frac{\partial Q_\lambda}{\partial p_{01}} \cdot \frac{\partial Q_\lambda}{\partial p_{23}} - \frac{\partial Q_\lambda}{\partial p_{02}} \cdot \frac{\partial Q_\lambda}{\partial p_{13}} + \frac{\partial Q_\lambda}{\partial p_{03}} \cdot \frac{\partial Q_\lambda}{\partial p_{12}}
=  t \cdot \left( p_{01}p_{23}-p_{02}p_{13} + p_{03}p_{12} \right)
\end{equation}
for some $t \in \mathbb{C}$.
This implies that $V(I)$ is isomorphic to the variety of all $Q \in \bbbp (\bbbc[\vec{p}]_2) \setminus \lbrace p_{01}p_{23}-p_{02}p_{13} + p_{03}p_{12} \rbrace$ satisfying~(\ref{eq:coisotropy_catanese}). Let $I_2$ be generated by the $2 \times 2$ minors of the $21 \times 2$ matrix that is obtained by deleting the middle column of the matrix in Figure~\ref{fig:bigmat}. Then $V(I_2)$ contains exactly those $Q \in \bbbp (\bbbc[\vec{p}]_2)$ satisfying~(\ref{eq:coisotropy_catanese}), and $V(I)$ is the projection of $V(I_2)$ 
from the center $(p_{01}p_{23}-p_{02}p_{13} + p_{03}p_{12})$. The ideal $I_2$ has codimension 11, degree 92, and is minimally generated by 20 quadrics. Interestingly, Catanese shows furthermore in \cite[Theorem 3.3]{catanese} that a hypersurface in $G(2,4)$ is coisotropic if and only if it is selfdual
in $\mathbb{P}^5$ with respect to the inner product given by the Pl\"ucker quadric.

\section{The Chow Variety}
\label{sec:chowvar}

In this section we study the Chow variety $G(2,2,4)$ of one-dimensional algebraic cycles of degree two in $\bbbp^3$. By~\cite[\S4.1, Ex.~1.3]{gkz}, the Chow variety $G(2,2,4)$ is the 
union of two irreducible components of dimension eight in $\bbbp^{19}$, one 
corresponding to planar quadrics and the other to pairs of lines. Formally, this means
that~$G(2,2,4) = V (P_{\mathrm{ChowConic}} ) \cup V (P_{\mathrm{ChowLines}})$,
where $P_{\mathrm{ChowConic}}$ is the homogeneous prime ideal in $\bbbq[V]$ whose variety comprises the Chow forms of irreducible curves of degree two in $\bbbp^3$. The ideal $P_{\mathrm{ChowConic}}$ has codimension 11 and degree 92, and it is minimally generated by 21 quadrics and 35 cubics.
The radical ideal $P_{\mathrm{ChowConic}} \cap P_{\mathrm{ChowLines}}$ has codimension 11, degree $232=92+140$, and  it is minimally generated by 230 cubics.

Since $G(2,2,4)$ should be contained in the coisotropic variety
$V(I)$, it seems that $P_{\mathrm{ChowConic}}$ is missing from the decomposition (\ref{coisotropy_ideal_decomposition}).  Here is the explanation:

\begin{proposition} \label{prop:chowhurw}
Every Chow form of a plane conic in $\bbbp^3$ is also a Hurwitz form.
In symbols,
$\,P_{\mathrm{Hurwitz}} \,\subset \,P_{\mathrm{ChowConic}} \,$ and thus 
$\,V ( P_{\mathrm{ChowConic}} ) \,\subset V ( P_{\mathrm{Hurwitz}}) $.
\end{proposition}

Our first proof is by computer: just check the inclusion of ideals in Macaulay2. 
For a conceptual proof, we consider a $4 \times 4$-symmetric matrix
 $M = M_0  + \epsilon M_1$, where ${\rm rank}(M_0) =1$.
By \cite[eqn.~(1)]{Stu}, the Hurwitz form
 of the corresponding quadric surface in $\bbbp^3\,$ is 
 $\,Q(\vec{p}) = \vec{p} (\wedge_2 M) \vec{p}^T$.
 Divide by $\epsilon$ and  let $\epsilon \rightarrow 0$.
 The limit is  the Chow form of the plane
  conic defined by restricting $M_1$ to ${\rm ker}(M_0) \simeq \bbbp^2$.
  This type of degeneration is familiar from the study of complete quadrics
  \cite{CGMP}.    Proposition \ref{prop:chowhurw}
  explains why the locus of irreducible curves
  is  not visible in~(\ref{coisotropy_ideal_decomposition}).

Gel'fand, Kapranov and Zelevinsky \cite[\S 4.3]{gkz} introduce
a class of differential forms in order to  discriminate Chow forms
among all coisotropic hypersurfaces.
In their setup, these forms represent the integrability of the $\alpha$-distribution
$\mathcal{E}_{\alpha,Z}$.  We shall apply the tools of computational
commutative algebra to shed some light on the
characterization of Chow forms via integrability of $\alpha$-distributions.

For this, we use local affine coordinates instead of Pl\"ucker coordinates. 
A point in the Grassmannian $G(2,4)$ is represented as the row space of
the matrix
\begin{equation}
\label{eq:twobyfour}
\left( \begin{array}{cccc}
1 &\, 0 &\, a_2 & a_3 \\
0 &\, 1 & \, b_2 & b_3
\end{array} \right) \; .
\end{equation}
We express the quadrics $Q$ in (\ref{general_quadric}) in terms of the local coordinates 
$a_2,a_3,b_2,b_3$, by substituting the Pl\"ucker coordinates with the minors of the 
matrix (\ref{eq:twobyfour}), i.e.,
\begin{equation}
p_{01} = 1, \,\,
p_{02} = b_2, \,\,
p_{03} = b_3, \,\,
p_{12} = -a_2, \,\,
p_{13} = -a_3, \,\,
p_{23} = a_2 b_3 - b_2 a_3 \; .
\end{equation}
We consider the following differential $1$-forms on affine $4$-space:
\[
\alpha_1^1 := \frac{\partial Q}{\partial a_2} d a_2 + \frac{\partial Q}{\partial a_3} d a_3,\quad
 \alpha_2^1 := \frac{\partial Q}{\partial a_2} d b_2 + \frac{\partial Q}{\partial a_3} d b_3,
\]
\[
\alpha_1^2 := \frac{\partial Q}{\partial b_2} d a_2 + \frac{\partial Q}{\partial b_3} d a_3,\quad
 \alpha_2^2 := \frac{\partial Q}{\partial b_2} d b_2 + \frac{\partial Q}{\partial b_3} d b_3 \; .
\]
By taking wedge products, we derive the $16$ differential $4$-forms
\begin{equation}
d Q \wedge d \alpha^i_j \wedge \alpha^k_l
\,\,=\,\, q_{ijkl} \cdot d a_2 \wedge d a_3 \wedge d b_2 \wedge d b_3
\,\,\quad \hbox{for $i,j,k,l \in \lbrace 1,2 \rbrace$} \; .
\end{equation}
Here the expressions $q_{ijkl} $ 
are certain polynomials in $ \bbbq[V][a_2,a_3,b_2,b_3]$.

Theorems~3.19 and~3.22 in~\cite[\S4.3]{gkz} state that a squarefree coisotropic
quadric $Q$ is a Chow form if and only if all~16 coefficients~$q_{ijkl}$ are
multiples of~$Q$.  By taking normal forms of the polynomials $q_{ijkl}$ modulo
the principal ideal $\langle Q \rangle$, we obtain a collection of $720$
homogeneous polynomials in $\vec{c}$.  Among these,~$58$ have degree three,~$340$
have degree four, and~$322$ have degree five.  The aforementioned result
implies that these~$720$ polynomials cut out~$G(2,2,4)$ as a subset of~$\bbbp^{19}$.

The {\em integrability ideal} $J \subset  \bbbq[V]$
is generated by these $720$ polynomials and their analogues from
other affine charts of the Grassmannian, obtained by permuting 
columns in~(\ref{eq:twobyfour}).  We know that $V(J)$ equals the union of $G(2,2,4)$ with all double hyperplanes in $G(2,4)$ (corresponding to $P_{\mathrm{Squares}}$)
set-theoretically.   Maple, Macaulay2 and Magma verified for us that it holds
scheme-theoretically:

     \begin{proposition}
The integrability ideal~$J$ is minimally generated by~210 cubics.
Writing $\mathfrak m$ for the irrelevant ideal~$\langle c_0,c_1,\ldots,c_{20} \rangle$ 
of~$\bbbq[V]$, we have
\begin{equation}
\sqrt{J}\,\, =\,\, (J : \mathfrak m)\, \,=\,\, P_{\mathrm{ChowConic}} \,\cap\, P_{\mathrm{ChowLines}} \,\cap\, P_{\mathrm{Squares}} \; .
\end{equation}
   \end{proposition}  

\section*{Conclusion}

We reported on computational experiments with hypersurfaces in the Grassmannian
$G(2,4)$ that are associated to curves and surfaces in $\bbbp^3$.  For degree
$2$, all relevant parameter spaces were described by explicit polynomials in
$20$ variables. 
All ideals and computations discussed in this note can be obtained at
\begin{center}
  \url{www3.math.tu-berlin.de/algebra/static/pluecker/}
\end{center}

Many possibilities exist for future work.  Obvious next
milestones are the ideals for the Chow varieties of degree $3$ cycles in
$\bbbp^3$,  and degree $2$ cycles in~$\bbbp^4$.  Methods from  representation
theory promise  a compact  encoding of their generators, in terms of
irreducible ${\rm GL}(4)$-modules.  Another question we aim to pursue is
motivated by the geometry  of condition numbers \cite{condition}: express the
volume of a tubular neighborhood of a coisotropic quadric in   ${\rm G}(2,4)$
as a function   of~$\vec{c}$.

\bibliographystyle{splncs}

\end{document}